# EXPLICIT COMPUTATION OF GALOIS P-GROUPS UNRAMIFIED AT P

Nigel Boston and Charles Leedham-Green

## 0. Introduction.

In this paper we introduce a new method for finding Galois groups by computer. This is particularly effective in the case of Galois groups of $p$-extensions ramified at finitely many primes but unramified at the primes above $p$. Such Galois groups have been regarded as amongst the most mysterious objects in number theory [10]. Very little has hitherto been discovered regarding them despite their importance in studying $p$-adic Galois representations unramified at $p$. The conjectures of Fontaine-Mazur [7] say that they should have no $p$-adic analytic quotients (equivalently, the images of the Galois representations should always be finite), and there are generalizations due to the first author [3], [4], suggesting that they should instead have 'large' actions on certain trees.

Below, we indicate how the method works in the case of 2-extensions of $\mathbf{Q}$. This case is chosen since powers of 2 grow relatively slowly allowing us extensive computations. The idea is to modify the $p$-group generation algorithm [9] so that as it goes along, it uses number-theoretical information to eliminate groups that cannot arise as suitable quotients of the Galois group under investigation. In the best cases we obtain a short list of candidates for the Galois group. Two phenomena are worth noting. First, as noted in [5], the same groups turn up repeatedly. Second, even when the short list consists of more than one group, the candidates typically are so similar that most questions regarding them will yield the same answer. In this way, we can tell a lot about the Galois groups without always pinning them down exactly.

Two items of progress are worth noting here. First, our exhaustive search allows us to go well beyond what was known previously, which was in general only the maximal nilpotency class 2 quotient of these Galois groups. Second, we obtain at least a conjectural characterization of the Galois group of the maximal 2-extension unramified outside $\{p, q\}$ in most cases when $p \equiv 3 \pmod 4$ and $q \equiv 5 \pmod 8$.

## 1. Some Basics.

Let $S$ be a finite set of odd primes. The maximal 2-extension of $\mathbf{Q}$ unramified outside $S$ (allowing ramification at infinity) will be denoted by $\mathbf{Q}_S$ and its Galois

The first author was partially supported by NSF grant DMS 99-70184. The authors thank Joann Boston for her help with the figures.

group over $\mathbf{Q}$ by $G_S$. We seek to find $G_S$ given a set $S$. We begin with some elementary results.

**Theorem 1.** (1) For every finite index subgroup $H$ of $G_S$, $H/H'$ is finite.
(2) [8] The generator rank $d(G_S)$ equals the relator rank $r(G_S)$. Equivalently, $G_S$ has trivial Schur multiplicator.

*Proof.* (1) follows by class field theory, from the finiteness of certain ray class groups. (2) is a theorem of Shafarevich.

**Corollary.** (1) If $\#S = 1$, then $G_S$ is a finite cyclic group.
(2) If $\#S \geq 4$, then $G_S$ is infinite. (Golod-Shafarevich.)

In [5], the first author and Perry exploited various circumstances in which a decomposition subgroup has index 1 or 2 in $G_S$. This led to the following families. (Here $S = \{p, q\}$.)

**Theorem 2.** Assume that $p \equiv 3 \pmod 4$, $q \equiv 3 \pmod 4$, and without loss of generality (by quadratic reciprocity) $p$ is a quadratic residue of $q$. If $2^k$ is the highest power of 2 dividing $q^2 - 1$, then

$$G_S \cong \langle\, a,\ b \mid a^2 = b^{2^k} = 1, \quad b^a = b^{-1+2^{k-1}} \rangle,$$

the semidihedral group of order $2^{k+1}$. Further, we can take $a$ to be complex conjugation and $b$ to be the generator of any inertia subgroup at $q$ of $G_S$.

**Theorem 3.** Assume that $p \equiv 3 \pmod 4$, $q \equiv 1 \pmod 4$, and $\left(\frac{p}{q}\right) = -1$. If $2^k$ is the highest power of 2 dividing $q - 1$, then

$$G_S \cong M_{k+2}(2) = \langle\, a,\ b \mid a^2 = b^{2^{k+1}} = 1, \quad b^a = b^{1+2^k} \rangle,$$

the modular group of order $2^{k+2}$. Further, we can take $a$ to be complex conjugation and $b$ to be the generator of any inertia subgroup at $q$ of $G_S$.

**Theorem 4.** Assume $p \equiv 3 \pmod 4$, $2^k (k \geq 2)$ exactly divides $q-1$, and either
(1) $k = 2$ and $p \in ({\mathbf{F}_q}^\times)^4$,
or (2) $k \geq 3$ and $p \in ({\mathbf{F}_q}^\times)^2 - ({\mathbf{F}_q}^\times)^4$.
Then

$$G_S \ \cong\ P_k \ :=\ \langle\, a,\ b \mid a^2 \ =\ b^{-1}ababab^{2^k-1}a \ =\ 1 \,\rangle$$

(presented as a pro-2 group), of order $2^{3k+1}$. Further, we can take $a$ to be complex conjugation and $b$ to be the generator of any inertia subgroup at $q$ of $G_S$.

Note the tendency for $G_S$ to belong to a small family of possibilities. This paper arose from the desire to study cases where $G_S$ does not have a decomposition subgroup of small index. (It can be checked that the groups given below do not have any metacyclic subgroups of small index. Note that decomposition subgroups must be metacyclic here since all ramification is tame.)

## 2. The Method.

The idea is to modify O'Brien's $p$-group generation algorithm [9], which allows computer algebra systems such as MAGMA [2] to find systematically all $p$-groups with some property, e.g. all 2-groups with 2 generators and order $\leq 1024$. Our modification is to have the program save those descendants that satisfy various conditions coming out of the number theory.

O'Brien's approach is to define, given $p$-group $G$, a sequence of characteristic subgroups by $P_1(G) = G,\ P_i(G) = P_{i-1}(G)^p(G, P_{i-1}(G))$. Thus, $G = P_1(G) \geq P_2(G) \geq ...$ The smallest $c$ such that $P_{c+1}(G) = \{1\}$ is called the *p-class* of $G$. Since the $i$-th term $\gamma_i(G)$ of the lower central series of a $p$-group $G$ is contained in the $i$-th term $P_i(G)$ of the lower $p$-central series, the $p$-class is bounded below by the nilpotency class. For most $p$-groups in question, $P_i(G) = \gamma_i(G)$ for all but a few initial values of $i$, and so the $p$-class and nilpotency class are equal.

A $p$-group $H$ is called a *descendant* of $G$ if $H/P_{c+1}(H) \cong G$, where $c$ is the $p$-class of $G$. It is an *immediate* descendant if it has $p$-class $c+1$.

O'Brien's algorithm finds all immediate descendants of a given $p$-group $G$. These can be conveniently arranged in the form of a tree. For instance, the immediate descendants of $C_2 \times C_2$ look like the following.

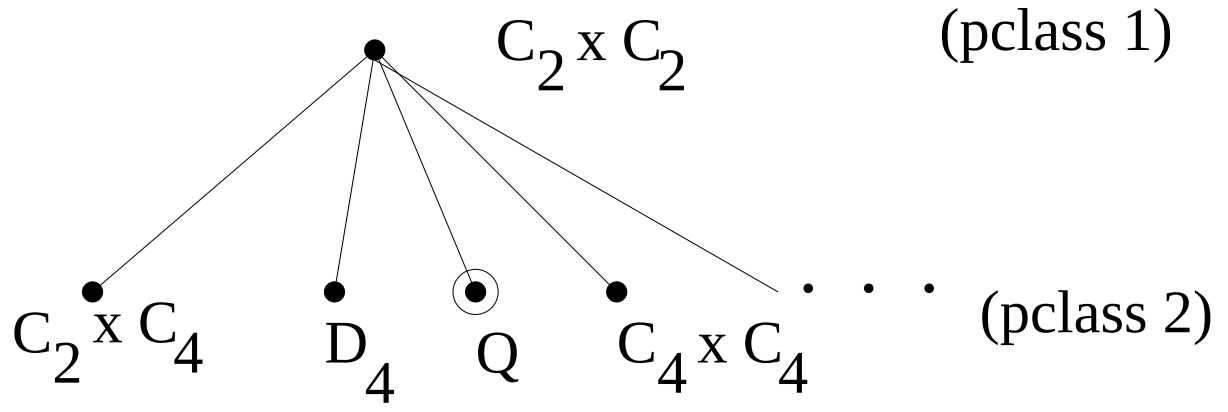

This can then be iterated to give all $p$-class 3, then $p$-class 4, etc. descendants of $C_2 \times C_2$. The tree contains each 2-group with 2 generators exactly once (up to isomorphism). The circle by the quaternion group, $Q_8$, indicates that it is 'terminal', i.e. has no descendants. The infinite ends of the tree yield infinite pro-2 groups; for instance $C_2 \times C_2 \leftarrow D_4 \leftarrow D_8 \leftarrow ...$ has inverse limit the dihedral pro-2 group.

Our aim then is, given $S, c$, to compute a short list (hopefully just one group) of possibilities for the $p$-class $c$ quotient $G_S/P_{c+1}(G_S)$.

We impose three kinds of number-theoretical constraint on $G_S$.

(I) (Top-end information) We can compute $H/H'$ for $H$ of small index in $G_S$ (index $\leq 2$ by class field theory, $\leq 16$ by machine, e.g. KASH [6]).

(II) (Local information) We know that $G_S$ has inertial generators $\{\tau_p : p \in S\}$ with the property that $\tau_p$ is conjugate to $\tau_p^p$. It also contains a complex conjugation $\tau_\infty$ of order 2.

(III) (Multiplicator information) $r(G_S) = d(G_S)$.

Note that if we had been considering $p$-extensions ramified at $p$, then $H/H'$ need not have been finite, which would make our approach much less effective.

## 3. An Example.

The first case not covered by the explicit families above is that of $S = \{5, 19\}$. We attempt to compute $G_S$ by O'Brien's program using the number-theoretical constraints above.

First, we show that $G_S/P_4(G_S)$ is a certain group $P$ of order 64. This is done by showing that $x^{16} - 70x^{14} + 1161x^{12} + 5125x^{10} + 8910x^8 + 7783x^6 + 5235x^4 + 593x^2 + 25$ has $P$ as Galois group, that its root field is unramified outside $S$, and that $P$ has 2-class 3 and is not a quotient of any other 2-class 3 2-group.

Second, we find the conjugacy classes of $P$ where $\tau_5, \tau_{19}, \tau_\infty$ can lie.

Third, we use KASH to compute $H/H'$ for all subgroups $H$ of index $\leq 4$ in $G_S$.

All the above goes in an input file. Our program now finds all immediate descendants of $P$ and tests them for 3 properties:

(i) that they contain elements $x_5, x_{19}$ that are respectively conjugate to their 5th and 19th powers and $x_\infty$ of order 2 such that their images in $P$ are in the conjugacy classes found above,

(ii) that the abelianizations of their subgroups of index $\leq 4$ are smaller than (i.e. quotients of) the abelianizations of the corresponding subgroups of $G_S$,

(iii) that their nuclear rank differs from their 2-multiplicator rank by at most 2 (see [9] for definitions of these terms).

Each group that passes all three tests is then saved, its immediate descendants computed, and the tests repeated on these. The search is done by depth first so as to save on storage. For large enough 2-class (ii) can be refined, since at this depth the abelianizations must actually equal those of the corresponding subgroups of $G_S$.

This continues until we reach groups with abelianizations of their subgroups of index $\leq 4$ exactly right and with trivial Schur multiplicator (which implies that the group is terminal). These groups are then candidates for $G_S$.

If the process terminates, then we know that $G_S$ is finite and is one of the candidate groups. It might happen, however, that the process does not terminate in a reasonable length of time. In that case, it could be that there are infinite groups satisfying our number-theoretical constraints or that there are only finite groups but their $p$-class is so large that they have yet to be reached. The input file to our program contains a limit for the $p$-class, which can be adjusted. In the case of $\{5, 19\}$, the process terminates in groups of 2-class 11.

Here is what the tree for $S = \{5, 19\}$ looks like. The numbers next to the nodes give the exponents of the orders. A circle around a node means that none of the group's descendants pass all three tests and that the group itself cannot be $G_S$. A square around a node means that we have reached a group that satisfies everything $G_S$ is known to satisfy.

The Tree for S={5,19}.

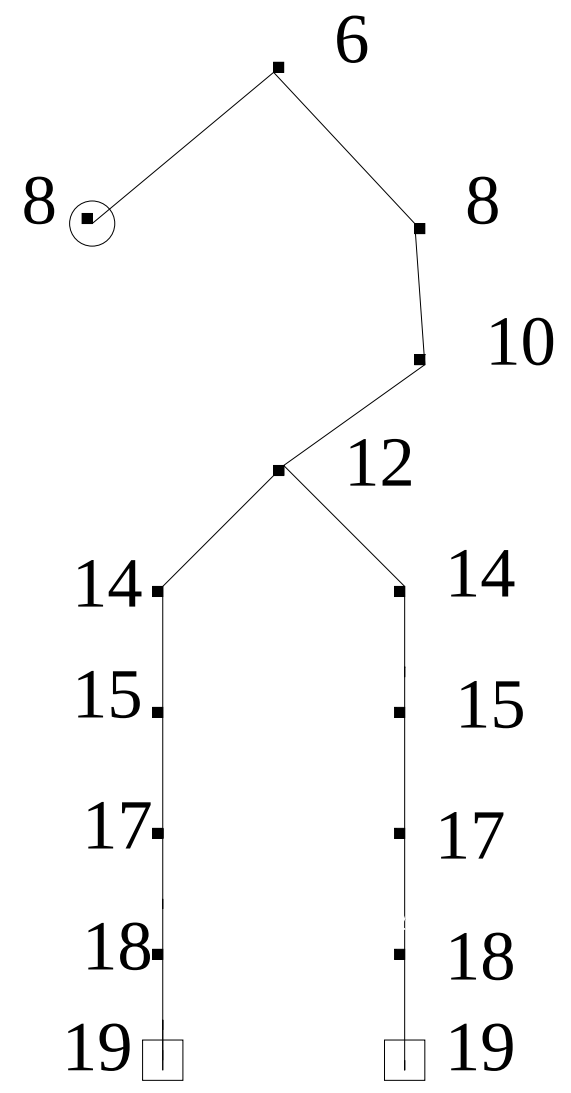

**Theorem 5.** If $S = \{5, 19\}$, then $G_S$ is of order $2^{19}$ and 2-class (indeed nilpotency class) 11. In fact, it is isomorphic to $\langle a, b : a^2, bab^2ab^{-5}ab^5ab^9ab^{-1}ab^5ab^{-4}a\rangle$ or to $\langle a, b : a^2, b^{-7}ab^{-6}ab^3ab^{-3}abab^{-1}ab^{-3}ab^{-4}a\rangle$.

These two groups are very similar to each other. For instance, their lattices of subgroups of index $\leq 32$ (even with abelianizations attached) match. This makes them indistinguishable by top-end methods, i.e. by calculation in low degree subfields of $\mathbf{Q}_S$. The advantage of this, though, is that the answer to many questions about $G_S$ can be given, since it is the same for either group.

**Remarks.** (1) The largest metacyclic subgroup of $G_S$ has order $2^9$ and so decomposition subgroups have index at least $2^{10}$. This means that the methods used in [5] for the earlier explicit families simply do not carry over.

(2) In fact, if we only use abelianizations for subgroups of index $\leq 2$ (which can be computed by class field theory rather than by machine), then we get the same tree and the same two groups output.

(3) Several other choices of $S$ lead to the same tree and indeed the same pair of groups, e.g. $\{3, 37\}, \{3, 61\}, \{5, 59\}, \{5, 139\}, \{11, 53\}, \ldots$.

Next we compute some other trees. For $S = \{13, 23\}$, we get the following. The same tree and groups are output for $S = \{13, 103\}, \{5, 199\}, \ldots$.

The Tree for S={13,23}.

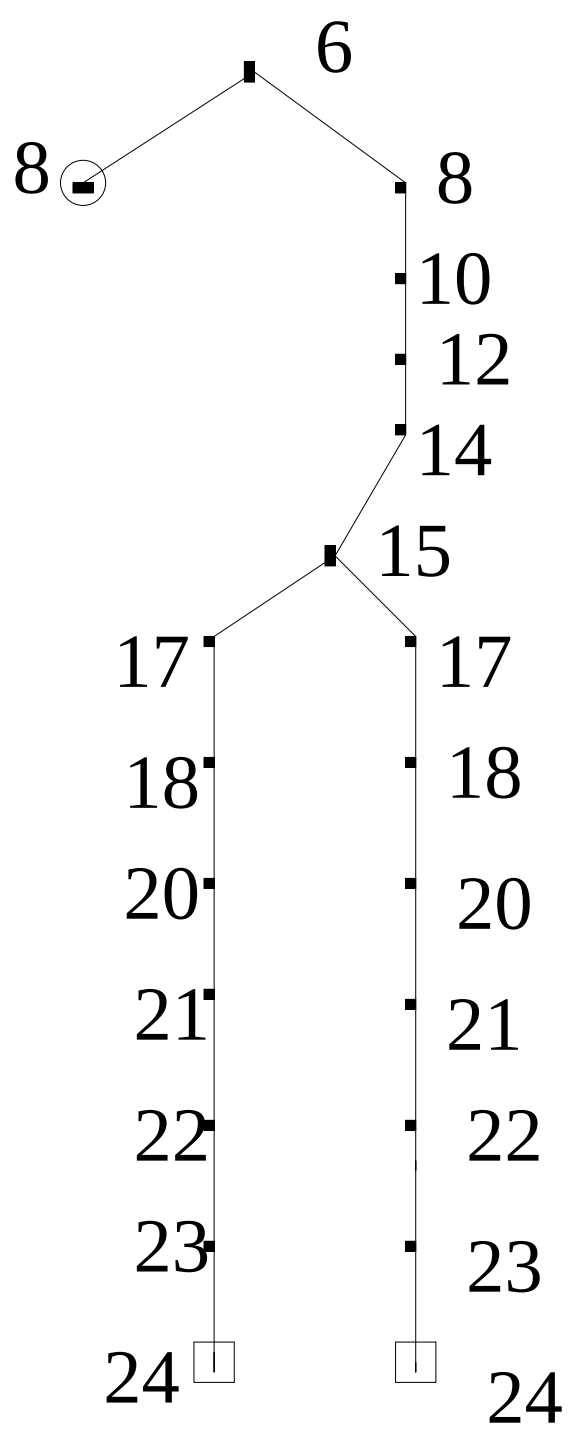

**Theorem 6.** If $S = \{13, 23\}, \{13, 103\}, \{5, 199\}, \ldots$, then $G_S$ is of order $2^{24}$ and 2-class (indeed nilpotency class) 15.

The tree for $S = \{5, 79\}$ is as follows.

The Tree for S={5,79}.

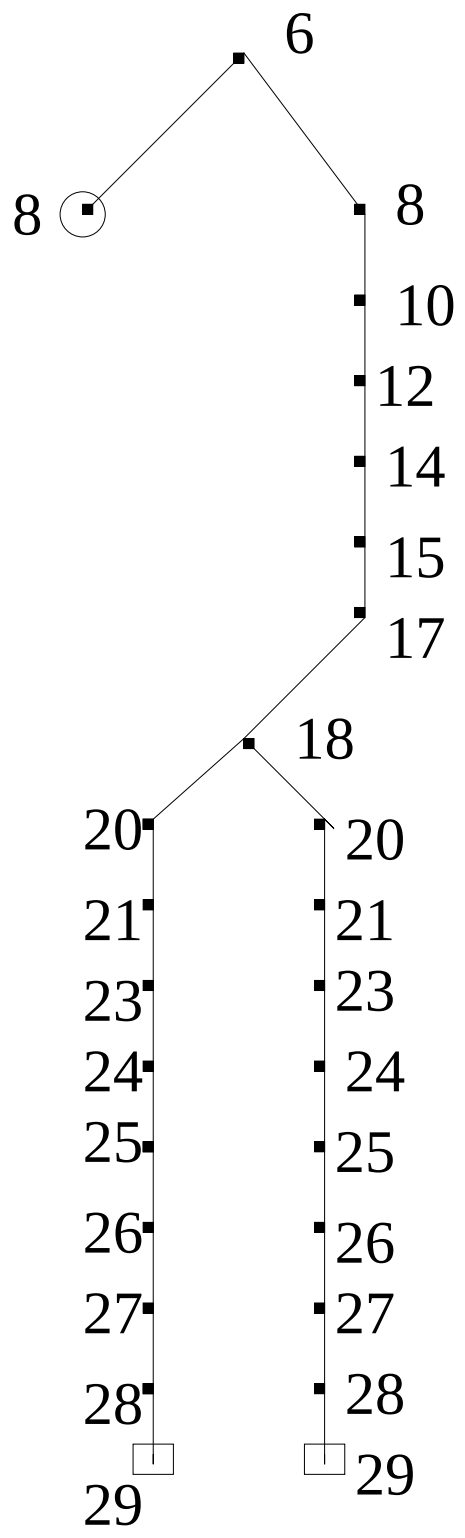

**Theorem 7.** If $S = \{5, 79\}$, then $G_S$ is of order $2^{29}$ and 2-class (indeed nilpotency class) 19.

It should be noted that no short presentations of these groups have been found. These trees and their output groups seem to fit into a pattern that we now formalize as a conjecture. (This pattern has been checked experimentally as far as $S = \{13, 127\}$, i.e. $k = 7$.)

First, note that by class field theory, if $p \equiv 3 \pmod 4$, $q \equiv 5 \pmod 8$, and $p$ is a square but not a 4th power modulo $q$, then the abelianizations of the subgroups of index 2 are $[4, 4], [2, 2, 2], [2, 2^n](n \geq 4)$, corresponding to $\mathbf{Q}(\sqrt{-p}), \mathbf{Q}(\sqrt{q}), \mathbf{Q}(\sqrt{-pq})$ respectively.

Suppose that $n = 4$ (as it does in all the cases above).

**Conjecture.** There is a family of groups $G(k)$ $(k \geq 2)$ such that $G(k)$ has order $2^{5k+9}$ and nilpotency class $4k + 3$, and such that if $2^k$ exactly divides $p + 1$ (together with the other assumptions above), then $G_S \cong G(k)$.

The main problem is to identify the groups $G(k)$. There are some similarities such as those noted above. Also, the number of generators of their normal subgroups apparently never exceeds 4 and there is always a unique maximal abelian normal subgroup $A(k)$ of order $2^{3k+3}$. On the other hand, the structures of $G(k)'/G(k)''$ and $A(k)$ vary erratically. For $k = 2, 3, 4, 5$, $G(k)'/G(k)''$ is respectively $[2, 4, 16], [8, 8, 8]$, $[4, 16, 32], [4, 32, 64]$, whereas $A(k)$ is respectively

[4, 4, 8, 32], [4, 8, 32, 32], [4, 16, 64, 64], [4, 16, 256, 256].

## 4. The Case $S = \{p, q, r\}$.

We now produce some cases with $\#S = 3$ and $G_S$ infinite. The Fontaine-Mazur conjecture and its generalizations apply to this case. Our method yields large quotients of $G_S$ against which to test out these conjectures.

Suppose $p, q, r \equiv 3 \pmod 4$. This keeps $G_S/G_S'$ small, namely $C_2 \times C_2 \times C_2$. The theory now comes in two flavors according as

(a) $\{p, q, r\}$ can be ordered so that $(\frac{p}{q}) = (\frac{q}{r}) = (\frac{r}{p})$ or

(b) they cannot.

### Case (a).

In this case, $G_S$ is apparently always finite. Consider, for example, $S = \{3, 11, 19\}$. The abelianizations for the seven subgroups of index 2 consist of [4, 4] once, [2, 2, 4] three times, and [2, 2, 8] three times. Taking as input the 2-class 1 quotient $P = G_S/P_2(G_S) \cong C_2 \times C_2 \times C_2$, our method yields the following tree.

Tree for S={3,11,19}.

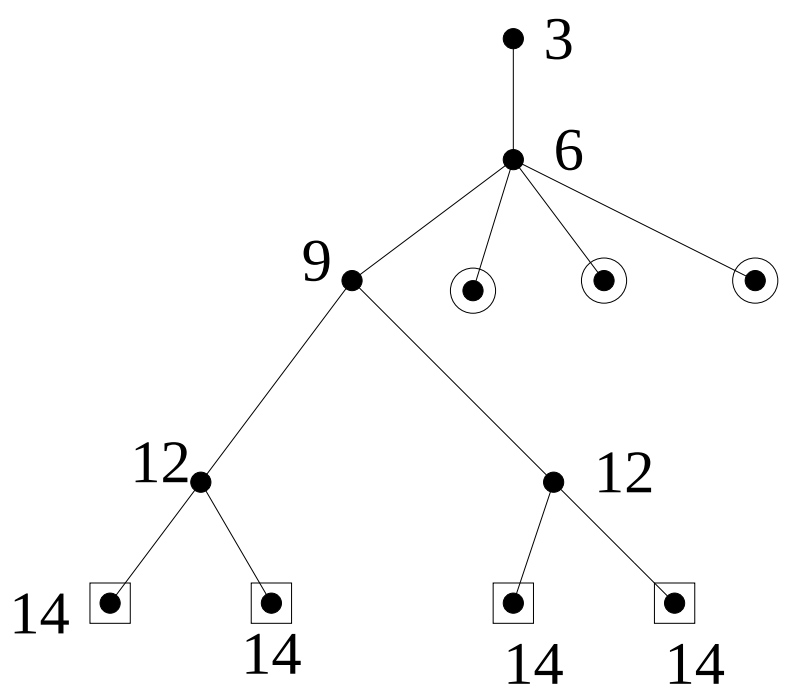

**Theorem 8.** If $S = \{3, 11, 19\}$, then $G_S$ is of order $2^{14}$ and 2-class (indeed nilpotency class) 5.

Once again, it is hard to tell which of the four groups output is actually $G_S$, but once again it does not matter for many purposes since the answer to many questions about them yields the same answer.

### Case (b).

In this case, $G_S$ is always infinite. This is because a subgroup of index 2 has abelianization [2, 2, 2, 2] and has as fixed field an imaginary quadratic field and then we apply the result of Shafarevich that such a subgroup must have 4 generators and 4 relations and hence be infinite by Golod-Shafarevich.

Let e.g. $S = \{3, 19, 43\}$. In this case the abelianizations of the index 2

subgroups are as small as they can be, namely $[4,4]$ once, $[2,2,4]$ three times, and $[2,2,2,2]$ three times. If we start with the 2-class 2 quotient $P = G_S/P_3(G_S)$ of order 64, then we get a tree as follows:

Tree for S={3,19,43}.

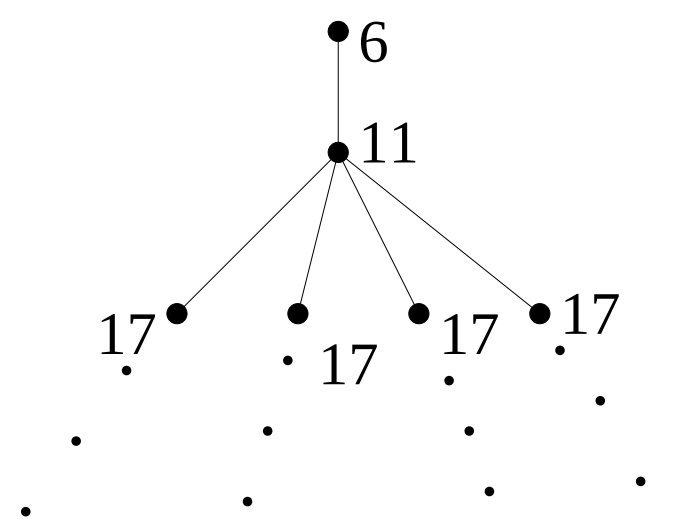

**Comments.**(1) The method has been used successfully in the case of 3-groups. For instance, in the case of three primes $S = \{p, q, r\}$ that are 1 (mod 3) but not 1 (mod 9), the first author has computed many cases of the Galois group of the maximal 3-extension unramified outside $S$. In each case, we obtain a finite group, in agreement with the result of Andoskii-Cvetkov [1], which says that this group is either finite or isomorphic to $\Gamma := \ker(SL(2,\mathbf{Z}_3) \to SL(2,\mathbf{F}_3))$. The latter would contradict the Fontaine-Mazur conjecture. This situation is explored in greater detail in joint work by the first author and Farshid Hajir.

(2) The Andoskii-Cvetkov groups have the interesting property that they form a family of 3-groups (in fact with 3 replaced by any odd prime) $G := G(n)$ such that $G/P_n(G)$ is isomorphic to $\Gamma/P_n(\Gamma)$ where $\Gamma$ is the above congruence subgroup. This suggests that in general the groups arising might fall into families with $G/P_n(G)$ isomorphic to $H/P_n(H)$ for some infinite 'governing' group $H$. For instance, for our family in theorem 4, such a role is played by $H = \langle a, b \mid a^2 = b^{-1}ababab^{-1}a = 1\rangle$ (we simply take the 2-adic limit of the given relation). This is the 2-adic space group $\mathbf{Z}_2 \wr C_2$. Unfortunately, there is no obvious candidate for such a governing group for our family in Section 3.

This unfortunately means that $\Gamma$ satisfies the three key properties we use to find our Galois groups. In other words, there exists $S$ such that the top-end information matches that of $\Gamma$, and it is easy to see that $\Gamma$ is generated by elements $x$ such that $x$ is conjugate to $x^n$ for some $n > 1$ and that $d(\Gamma) = r(\Gamma) = 3$. Perhaps there is some additional piece of information we have overlooked that excludes $\Gamma$ from being one of our Galois groups (in accord with the Fontaine-Mazur conjecture).

(3) As Hajir has noted, $G_S$ can be infinite even if $\#S = 2$, for instance, if $S = \{17, 101\}$. This follows by applying Golod-Shafarevich to the maximal real subfield of $\mathbf{Q}(\zeta_{17})$. This introduces the interesting possibility of constructing a family of examples for which the first few are finite and can be found explicitly but the later ones are infinite. For instance, $\{17, 101\}$ might be in a family with general term $\{p, 6p-1\}$ with $2^k$ exactly dividing $p-1$. The first term is e.g. $S = \{5, 29\}$ and we intend to compute $G_S$ once some computer glitches are overcome. Replacing the 4 by 16 (somehow) in the presentation yields a guess for $G_S$ when $S = \{17, 101\}$.

## 5. Impact on the Fontaine-Mazur Conjecture and its Generalizations.

Let $W$ be the automorphism group of the binary rooted tree. In other words, we set $W_1 = C_2$ and $W_n = W_{n-1} \wr C_2$, so that $W_n$ is isomorphic to the Sylow 2-subgroup of the symmetric group on $2^n$ letters. The exponent of its order is $2^n - 1$. Let $W = \lim_{\leftarrow} W_n$.

If $G$ is a closed subgroup of $W$, its Hausdorff dimension is defined to be $\liminf_{n\to\infty} \frac{\log|G_n|}{\log|W_n|}$, where $G_n$ is the image of $G$ in $W_n$. The first author has made the following conjecture regarding $G_S$ [4]. (Note that a just-infinite pro-$p$ group is one that is infinite but all of its proper quotients are finite.)

**Conjecture.** The just-infinite quotients of $G_S$ are isomorphic to closed subgroups of $W$ of nonzero Hausdorff dimension.

Note that this implies the Fontaine-Mazur conjecture in this case, in that a 2-adic analytic group cannot embed in $W$ with nonzero Hausdorff dimension. (The conjecture amounts to saying that the just-infinite quotients of $G_S$ are Grigorchuk's branch groups, together with a conjectural characterization of branch groups in terms of Hausdorff dimension.)

The conjecture also implies that if $G_S$ is infinite, then there exists a representation $G_S \to W$ with image of nonzero Hausdorff dimension. We use our new computational tools to test this out.

**Example.** Let $S = \{3, 19, 43\}$. Let $f = x^{16} - 4x^{15} + 10x^{13} + 14x^{12} - 42x^{11} - 20x^{10} + 80x^9 + 15x^8 - 100x^7 + 10x^6 + 78x^5 - 25x^4 - 26x^3 + 9x^2 - 10x + 13$. Its Galois group is a 2-group since its root field has subfields of every index. The splitting of the polynomial modulo various primes tells us by van der Waerden's method that its Galois group contains elements with various cycle structures. Using the MAGMA database of transitive groups of degree 16, we find there is exactly one subgroup of $W_4$ that contains elements with all these cycle structures and that has abelianization $C_2 \times C_2 \times C_2$. This then gives us a map from $G_S \to W_4$ with image TransitiveGroup$(16, 1735)$ of order $2^{13}$.

Department of Mathematics, University of Illinois, Urbana, Illinois 61801
*E-mail address*: boston@math.uiuc.edu

School of Mathematical Sciences, Queen Mary and Westfield College, University of London, Mile End Road, London E1 4NS
*E-mail address*: C.R.Leedham-Green@qmw.ac.uk